\documentclass[12pt,english]{article}
\def\IC{\relax\,\hbox{$\inbar\kern-.3em{\rm C}$}}
\newcounter{lemma}[section]
\newcounter{proof}[section]
\newcounter{claim}[section]
\newcounter{corollary}[section]
\newcounter{theorem}[section]
\newcounter{proposition}[section]
\newcounter{definition}[section]
\newcounter{example}[section]
\newcounter{problem}[section]
\newcounter{remark}[section]
\newcounter{question}[section]
\def\vflx#1#2{\mathrel {\buildrel \hbox{$#1$} \over {\hbox to #2em
{\rightarrowfill}}}}
  \def\flx#1{\vflx{#1}{2}}
  \def\flxb{\Big\downarrow}

\date{}
\title{When Aut($\cal A$) and Homeo(Prim($\cal A$)) are homeomorphic, where $\cal A$ is a $C^*$-algebra
}
\author{B. BOUALI }

\begin{document}
\maketitle
%\def\thefootnote{\fnsymbol{footnote}}
%\footnotetext{}
\abstract{In this paper, we discuss when Aut($\cal A$) and Homeo(Prim($\cal A$)) are homeomorphic, where $\cal A$ is a $C^*$-algebra. 
}
\vspace{10mm}\\
A. M.S 2000 subject classification  : 46H05
\vspace{5mm}
\section{Preliminaries}
Let  $\cal A$ be $C^*$-algebra. Then the collection Aut($\cal A$) of $\star$-automorphisms of $\cal A$ is a group under composition. We give  Aut($\cal A$) the point-norm topology, that is, $\alpha_n \mapsto \alpha$ if and only if $\alpha_n (a) \mapsto \alpha (a) $ for all $a\in \cal A$. Then Aut($\cal A$) is a topological group.

Let $(X,{\cal X})$ and $(Y,{\cal Y})$ be topological spaces, and let $\phi : X \mapsto Y$ be a map. Then $\phi$ is called continuous (from $(X,{\cal X})$ to $(Y,{\cal Y})$) if the preimage $\phi^{-1} (V)$ of $V$ under $\phi$ belongs to $\cal X$ for each $V\in {\cal Y}$. A bijection between topological spaces is called a homeomorphism if it is continuous and has a continuous inverse.
\subsection{The hull kernel topology} The topology on Prim($\cal A$) is given by means of a closure operation. Given any subset $W$ of Prim($\cal A$), the closure $\overline{W}$ of $W$ is by definition the set of all elements in Prim($\cal A$) containing $\cap W=\lbrace \cap I  ~: ~ I \in W\rbrace $, namely 
$$
 \overline{W}=\lbrace I\in Prim({\cal A}): I \supseteq \cap W \rbrace$$

It follows that the closure operation gefines a topology on Prim($\cal A$) which called Jacobson topology or hull kernel topology [1].

\begin{proposition} The space Prim($\cal A$) is a $T_0$-space, i.e. for any two distinct points of the space there is an open neighborhood of one of the points which does not contain the other.\end{proposition}

\begin{proposition} If $\cal A$ is $C^{*}$-algebra, then Prim($\cal A$) is locally compact.If $\cal A$ has a unit, then Prim($\cal A$) is locally compact. \end{proposition}

\begin{remark} The set of ${\cal K}({\cal H})$ of all compact operators on the Hilbert space $\cal H$ is the largest two sided ideal in the $C^*$-algebra  ${\cal B}({\cal H})$ of all bounded operators.\end{remark}

\begin{definition}
A $C^*$-algebra $\cal A$ is said to be liminal if for every irreducible representation $(\pi, {\cal H})$ of $\cal A$ one has that $\pi(\cal A)= {\cal K}({\cal H})$
\end{definition}

So, the algebra $\cal A$ is liminal if it is mapped to the algebra of compact operators under any irreducible representation. Furthermore if $\cal A$ is a liminal algebra, then one can  prove that each primitive ideal of $\cal A$ is automatically a maximal closed two-sided ideal. As a consequence, all points of Prim($\cal A$) are closed and Prim($\cal A$) is a $T_1$-space . In particular, every commutative  $C^*$-algebra is liminal [1].
    
\section{Aut($A$) and Homeo(Prim($\cal A$))}
\begin{theorem}
Let $\cal A$ be liminary $C^*$-algebra, If $\alpha\in Aut(\cal A)$, there is $h \in Homeo(Prim(\cal A))$ such that :
$$
\alpha(a)(\pi)=h^{-1}(\pi)(a) ~~\forall a \in A \mbox{~and~} Ker(\pi) \in Prim(\cal A)$$
The map $\alpha \mapsto h$ is a homeomorphism. 
\end{theorem}

{\bf Proof} The primitive ideal of $\cal A$ is a maximal ideal [1 ,corollary 4.1.11.(ii)], Let $I_\pi$ be a maximal ideal of $\cal A$ for some $\pi \in {\cal A}$. Since $\alpha^{-1}(I_\pi)$ is a maximal ideal, there is a function h
$$ h: Prim({\cal A}) \mapsto Prim({\cal A})$$\\ such that $\alpha^{-1}(I_\pi)=I_{h(\pi)}$. 

Since we can replace $\alpha^{-1}$ by $\alpha$, it follows that $h$ is a bijection. we have induced isomorphism $\chi_\pi : {\cal A}/I_\pi \mapsto \bf C$ given by $\chi_\pi (a)=\pi(a)$  and  $\beta : {\cal A}/I_\pi \mapsto {\cal A}/I_{h^{-1}(\pi)}$ given by $\beta(a+I_\pi ) =\alpha (a) +I_{h^{-1}(\pi)}$. Therefore we get a commutative diagram:

$$
  \matrix{
  {\cal A}/I_{\pi} & \flx\beta & {\cal A}/I_{h(\pi)}\cr\cr
  \chi_\pi \flxb && \flxb \chi_{h(\pi)} \cr\cr
  {\bf C} & \flx{\gamma} & {\bf C}
  }
  $$
\vspace{5mm}\\
 and an induced  isomorphism $\gamma : {\bf C} \mapsto {\bf C}$ defined by $$\gamma (\pi(a))=h^{-1}(\pi)(\alpha(a))$$ Since $\gamma$ must be the identity maps , then  $\pi(a)=h^{-1}(\pi)(\alpha(a))$.

All opens set of Prim($\cal A$) are of the forms $U_I =\lbrace P\in Prim({\cal A}) : P   \not \supseteq  I\rbrace$
[2,3], calculate $h^{-1}(U_I)$.\\
$h^{-1}(U_I)= \lbrace  ~~~~\pi~~~  : ker(\pi) \in Prim({\cal A}) \mbox{~~and~~} ker(\pi)   \not \supseteq  I\rbrace$\\
$~~~~~~~~~~~~= \lbrace  h^{-1}(\pi) : ker(\pi)\in Prim({\cal A}) \mbox{~~and~~} ker(\pi)   \not \supseteq  I\rbrace$\\
$~~~~~~~~~~~~= \lbrace  ~~~~\pi '~~~ : ker(\pi ')\in Prim({\cal A}) \mbox{~~and~~} ker( h(\pi'))   \not \supseteq  I\rbrace$\\
$~~~~~~~~~~~~= \lbrace  ~~~~\pi '~~~ : ker(\pi ')\in Prim({\cal A}) \mbox{~~and~~} \alpha^{-1}(I_\pi)   \not \supseteq  I\rbrace$\\
$~~~~~~~~~~~~= \lbrace  ~~~~\pi '~~~ : ker(\pi ')\in Prim({\cal A}) \mbox{~~and~~} ker( \pi')   \not \supseteq  \alpha(I)\rbrace$\\$~~~~~~~~~~~~= U_{\alpha(I)}$.\\

Then $h^{-1}(U_I)$ is an open set and h is continuous. Replace $\alpha$ by $\alpha^{-1}$, it follows that $h^{-1}$ is continuous, then $h$ is a homeomorphism.\\
We still have to check that $\alpha \mapsto h$ is a homeomorphism. Let $\alpha_n $ in Aut($\cal A$) be given by  $\pi(\alpha_n(a))=h_n(\pi)(a)$.

Suppose that $\alpha_n  \mapsto  \alpha$. If, contrary there exists $\pi \in Homeo(Prim(\cal A))$ such that $h_n(\pi) \not \mapsto h(\pi)$  and $h_n(\pi)(a) \not \mapsto h(\pi)(a)$, then $\pi(\alpha_n(a)) \not\mapsto \pi(\alpha(a))$. Thus $\alpha_n(a)  \mapsto  \alpha(a)$ in $\cal A$. This is a contradiction, and we must have $h_n \mapsto h$ in Homeo(Prim($\cal A$). Since $\alpha_n^{-1}  \mapsto  \alpha^{-1}$, we have also $h_n^{-1}  \mapsto  h^{-1}$. Thus  
$h_n  \mapsto  h$ in Homeo(Prim($\cal A$)).

Now suppose that $h_n  \mapsto  h$ in Homeo(Prim($\cal A$)). We need to show that given $a\in \cal A$, and $I_\pi$ in Prim($\cal A$). If not, there exists $a\in \cal A$, such that $\alpha_n(a)  \not\mapsto  \alpha(a)$. Thus   $h_n(\pi)(a) \not \mapsto h(\pi)(a)$. This is a contradiction, and we must have $\alpha_n \mapsto \alpha$ in Aut($\cal A$).

\end{document}